# High-fidelity quantitative differential phase contrast deconvolution using dark-field sparse prior

Shuhe Zhang,[1,2,*] Tao Peng,[1] Zeyu Ke,[1] Meng Shao,[1] Tos T. J. M. Berendschot,[2] and Jinhua Zhou[1,3,**]

1. School of Biomedical Engineering, Anhui Medical University, Hefei 230032, China.
2. University Eye Clinic Maastricht, Maastricht University Medical Center + , P.O. Box 5800, Maastricht, 6202 AZ, the Netherlands
3. Anhui Provincial Institute of Translational Medicine, Anhui Medical University, Hefei 230032, China.

*shuhe.zhang@maastrichtuniversity.nl;
**zhoujinhua@ahmu.edu.cn

**Abstract:** Differential phase contrast (DPC) imaging plays an important role in the family of quantitative phase measurement. However, the reconstruction algorithm for quantitative DPC (qDPC) imaging is not yet optimized, as it does not incorporate the inborn properties of qDPC imaging. In this research, we propose a simple but effective image prior, the dark-field sparse prior (DSP), to facilitate the phase reconstruction quality for all DPC-based phase reconstruction algorithms. The DSP is based on the key observation that most pixel values for an idea differential phase contrast image are zeros since the subtraction of two images under anti-symmetric illumination cancels all background components. With this DSP prior, we formed a new cost function in which $L_0$-norm was used to represent the DSP. Further, we developed two different algorithms based on (1) the Half Quadratic Splitting, and (2) the Richardson-Lucy deconvolution to solve this NP-hard $L_0$-norm problem. We tested our new model on both simulated and experimental data and compare against state-of-the-art methods including $L_2$-norm and total variation regularizations. Results show that our proposed model is superior in terms of phase reconstruction quality and implementation efficiency, in which it significantly increases the experimental robustness, while maintaining the data fidelity.

## 1. Introduction

Differential phase-contrast microscopy (DPC), a non-interferometric quantitative phase retrieval approach, has been used for label-free and stain-free optical imaging of live biological specimens both *in vitro* [1-6] and *in vivo* [7, 8]. A quantitative DPC experimental layout involves a 4-f microscopy system in which a programmable LED or LCD illumination source generates anti-symmetric illumination patterns [6, 9]. With the combination of oblique plane wave illumination and low-pass filtering of the objective lens, the DPC converts the unmeasurable sample phase into a phase-contrast intensity image. By collecting at least 4 phase-contrast images with asymmetric illuminations in opposite directions, the phase component of the sample can be reconstructed through a non-blind deconvolution process, where the convolution kernel, in an ideal condition, is defined by the Fourier transform of phase contrast transfer function (PTF) [6]. The spatial deconvolution is then transformed into Fourier space division.

However, since no optical system is perfect, DPC raw images are always corrupt by noise, illumination fluctuations, and optical aberrations, resulting in the degeneration of the phase recovery results. The most comment problem in the DPC deconvolution process (and for all deconvolution processes) is that the PTF contains many pixels that approach zero (along its axis of asymmetry and beyond the passband) that will largely amplify the noise signal during the deconvolution process [6]. Therefore, a small constant is added to the denominator resulting in the Tikhonov-regularization, also known as $L_2$-norm regularization [6, 10-12]. This is widely used for almost all state-of-the-art DPC experiments. Besides, many methods have been proposed to improve the performance of DPC for both reconstruction quality and time resolutions, which can be categorized into two families: (1) PTF engineering, and (2) algorithm modification.

The DPC reconstruction is directly related to the PTF which is determined by the illumination pattern and pupil of the objective lens. Since it is inconvenient to directly modify the pupil of the objective lens, the PTF engineering is converted into illumination pattern optimization [13]. In these studies, illumination is optimized using, for example, radially asymmetric patterns [14], gradient amplitude patterns [15], or ring-shaped sine patterns [11]. Using the modified illumination pattern, the PTF can be better defined and the impact of noise is suppressed. Color-coded DPC based on wavelength multiplexing is able to achieve single-shot high speed DPC imaging [12, 16, 17]. For the algorithm modification, only few papers have been published and it is known that the optimizing illumination becomes the most direct and effective approach to improve the image quality of DPC [11]. Nevertheless, the modification of illumination pattern will limit the application of DPC, and the actual illumination patterns in DPC experiments are not precisely fit the theoretical assumptions, since there is no perfect light source.

To tackle the impact of noise and illumination fluctuations, and so significantly improving the reconstruction quality of DPC deconvolution, in this paper, we propose a high fidelity DPC reconstruction algorithm in which the dark-field sparse prior (DSP) and Hessian penalty are embedded. The new DPC reconstruction algorithm is termed SH-DPC which is suitable for almost all DPC experiments and does not require any modifications to the optical devices. A demonstration of the SH-DPC phase reconstruction results is shown in Fig. 1.

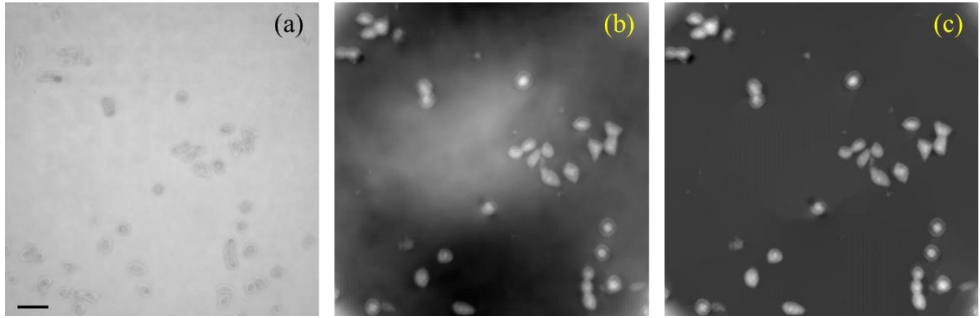

Fig. 1. DPC reconstruction using dark-field sparse prior. (a) bright-field image. (b) Reconstruction using TV-regularization. (c) Reconstruction using DSP. Scale bar is 65 μm.

In the following content, section 2 briefly review the off-the-shelf DPC reconstruction algorithms, and section 3 laid out the dark-field sparse prior which reveals the fact that the idea differential phase-contrast images are sparse and was hitherto neglected. In section 4, we use $L_0$-norm to quantify the image sparsity and formulate a new cost function for the DPC inverse problem, and two solutions for the $L_0$-norm problem are developed. Finally, section 5 shows the experimental results for validation and application of the SH-DPC and followed by discussion in section 6, and concluding remarks in section 7.

## 2. Related works

The inverse problem of DPC is an non-blind deconvolution task which is given by

$$\arg\min_{\varphi} \sum_{n=1}^{N} \left\| h_n \otimes \varphi - s_{n,dpc} \right\|_2^2 + \alpha \Phi(\varphi), \tag{1}$$

where $\varphi$ is the phase of the object to be solved. $h_n$ is the matrix denoting the $n$-th point spread function corresponding to the phase transfer function of $n$-th DPC image. $s_{n,dpc}$ is the $n$-th DPC image, which is calculated from the differences between two phase-contrast images under anti-systematic oblique illuminations

$$s_{n,dpc} = \frac{I_{n,r} - I_{n,l}}{I_{n,r} + I_{n,l}}, \qquad (2)$$

here, for example, $I_{n,r}$ and $I_{n,l}$ are captured images under anti-systematic oblique illuminations from right, and left, respectively. $\Phi(\varphi)$ in Eq. (1) is the penalty function on $\varphi$ to regularize the solution of Eq. (1), and $\alpha$ is the penalty parameter. Equation (1) is a non-blind deconvolution problem since $h_n$ are known according to the optical layouts.

For most of off-the-shelf algorithm $\Phi(\varphi) = \|\varphi\|_2^2$ is used to suppress the impact of noise which is corresponding to the $L_2$-norm/Tikhonov regularization, and $\alpha$ is manually adjusted according to the noise-level of the captured raw images [6, 11, 14, 16, 17, 19]. Ideally, $\alpha$ should be assigned a small value in order to keep the calculated $\varphi$ close to its ground truth value. However, when the raw images are corrupted by noise, $\alpha$ should be assigned a large value to restrict the noise signals, yielding the calculated $\varphi$ non-longer corresponding to its really value as the $L_2$-norm suppress the signals of both phase and noises simultaneously [20].

For better implement, one can uses the total-variant (TV) regularization (known also as gradient $L_1$-norm regularization) to suppress the noise signal where $\Phi(\varphi) = TV(\varphi) = \|\nabla\varphi\|_1$ [18, 21]. The TV regularization utilizes the first-order gradient of the image $\varphi$ and suppress the noise signals by erasing large values of gradient [21], however, may causes stair-case like artifact [22].

Although the above-mentioned methods are widely used for solving DPC inverse problem, they are not fully optimized according to specific feature of DPC raw data. In this study we propose the dark-field sparse prior and design an adaptive noise robust DPC solver by using the dark-field sparse and Hessian regularizations.

Our contribution of this work can be summarized as follows: (1) we present a new image prior termed as dark-field sparse prior (DSP); (2) We adopt $L_0$-norm on the DSP involved term, and provide an effective optimization scheme for the cost functions under two different frameworks; (3) our method performs well on both synthetic DPC data and real DPC experiments, against state-of-the-art algorithms. In the following section we introduce the new prior and prove why it works mathematically.

## 3. Dark-field sparse prior

The prior is based on a proposition of DPC data that (1) the DPC image is a sparse matrix, and (2) the sparsity of an idea noise-free DPC image is greater than noise corrupted DPC images (as shown in Fig. 2).

To better illustrate this observation, we consider the point spread function $h_n$ of the differential phase contrast image $s_{n,dpc}$, which is given by

$$h_n(\mathbf{x}) = \mathcal{F}^{-1}\left[H_{ph,n}(\mathbf{k})\right], \qquad (3)$$

where $\mathcal{F}^{-1}$ is the inverse Fourier transform. $\mathbf{k}$ is the coordinate vector in Fourier space. $H_{ph,n}(\mathbf{k})$ is the contrast phase transform function (PTF) and is determined by [6]

$$H_{ph,n}(k) = i \cdot \iint P(k_{ill})\left[S(k_{ill})P(k_{ill}+k) - S(k_{ill})P(k_{ill}-k)\right]d^2k_{ill}, \quad (4)$$

where $k_{ill}$ is the coordinate vector in illumination pupil plane. $S(k)$ is the illumination pupil and is a pure-real function. $P(k)$ for an idea objective lens is a pure-real even function. Let $k = -k$ in Eq. (4) we obtain $H_{ph,n}(k) = -H_{ph,n}(-k)$. Therefore, $H_{ph,n}$ is an odd function which meets:

$$H_{ph,n}(0) = 0. \quad (5)$$

Furthermore, according to the property of Fourier transform, $h_n(x)$ is a pure-real odd function thus we have

$$\int_{-\infty}^{\infty} h_n(x)dx = 0, \quad (6)$$

which implies the fact that the integration of an arbitrary odd function among a symmetric interval is zero.

The sparsity of $s_{n,dpc}$ can be analyzed from both Eq. (5) and Eq. (6). First, according to Eq. (5), the $H_{ph,n}$ completely removes the direct-current (DC) components in Fourier space which leads to the fact that the absolute value of DPC image has the similar visual performance to that of dark-field images (background is black). Second, according to Eq. (6), convoluting an image using $h_n(x)$ generates a sparse image as the kernel cancels out pixels that has the same values. In conclusion, $s_{n,dpc}$ in an ideal condition is a sparse matrix.

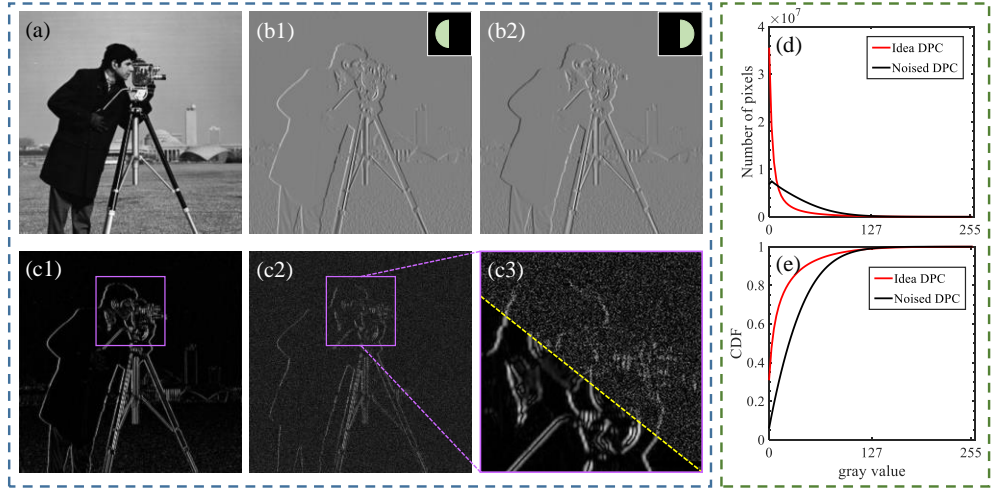

Fig. 2. Simulation of differential phase contrast images. (a) Phase pattern. (b1) and (b2) are simulated image under oblique illumination from left and right, respectively. (c1) is the absolute value of differential phase contrast image of (b1) and (b2). (c2) are the absolute value of (b1) and (b2) when corrupted by noises SNR = 0.7. (c3) enlarged part in the purple box. (d) and (e) are histogram and CDF (Cumulative Distribution Function) for 320 simulated DPC images with and without noise corruption.

The simulated DPC images are shown in Fig. 2. We generate the phase contrast images through simulation as shown in Fig. 2 (b1) and 2 (b2). Their difference pattern is shown in Fig. 2 (c1) in which the pattern is similar to a dark-field images as most of the pixels being zeros. If we add random noise to both Fig. 2 (b1) and 2 (b2), the difference pattern will be shown in Fig. 2 (c2). Partial enlarged image for the purple box is shown in Fig. 2 (c3). Since the noise is random and cannot be cancel by differential operation, the back grounds where they should be zero are now filled by noise pixels, and the difference pattern are not as sparse as the noise-free one.

We also design massive simulations based on different PTFs, and perform statistical analysis (please see **Supplementary Note 1** for details). As shown in Fig. 2(d) and 2(e), the idea DPC images are even sparser than noise corrupted DPC images. We therefore are inspired by this sparse observation and propose the dark-field sparse prior (DSP).

Mathematical proofs are as followed: taking $L_0$-norm on both side of Eq. (2) yielding

$$\left\| s_{n,dpc} \right\|_0 = \left\| \frac{I_{n,r} - I_{n,l}}{I_{n,r} + I_{n,l}} \right\|_0 = \left\| I_{n,r} - I_{n,l} \right\|_0. \tag{7}$$

The $L_0$-norm counts the nonzero elements of vector (matrix) and meets the triangle inequality in which $\left\| A + B \right\|_0 \leq \left\| A \right\|_0 + \left\| B \right\|_0$. Since the pixel values in $I_{n,r} + I_{n,l}$ are a non-negative, the $L_0$-norm of $s_{n,dpc}$ equals to the $L_0$-norm of $(I_{n,r} - I_{n,l})$. Let $\varepsilon_{n,r}$ and $\varepsilon_{n,l}$ be the noise signal or uneven background illumination imposed on $I_{n,r}$ and $I_{n,l}$. Since the random noises and uneven background cannot be completely cancel out by the subtraction, we are able to obtain the following inequation

$$\left\| I_{n,r} - I_{n,l} \right\|_0 \leq \left\| (I_{n,r} + \varepsilon_{n,r}) - (I_{n,l} + \varepsilon_{n,l}) \right\|_0 \leq \left\| s_{n,dpc} \right\|_0 + \left\| \varepsilon_{n,r} - \varepsilon_{n,l} \right\|_0. \tag{8}$$

Equation (8) implies that the present of noises, uneven background illumination, and illumination fluctuations tends to increase the $L_0$-norm of $s_{n,dpc}$ since those noises cannot be entirely canceled out during the subtraction. In other words, high quality DPC image tends to be sparser than those degraded one.

## 4. Proposed model: Sparse-Hessian DPC

In this section, we put forward a concrete DPC deconvolution model and two effective optimization scheme. Vectorizing $\varphi$ and $s_{n,dpc}$ within a conventional DPC framework, our sparse-Hessian cost function is defined as,

$$E(\varphi) = \sum_{n=1}^{N} \left\| K_n \varphi - s_{n,dpc} \right\|_2^2 + \alpha \sum_{n=1}^{N} \left\| K_n \varphi \right\|_0 + \beta \left\| \mathcal{H} \varphi \right\|_1, \tag{9}$$

where $K_n$ is the matrix denoting the convolution with kernel $h_n(x)$. The first fidelity term in the augmented Lagrangian Eq. (9) enforces similarity between convolution result $K_n \varphi$ and the observed DPC images $s_{n,dpc}$. $\alpha$ and $\beta$ are positive penalty parameters for the following regularization terms and are determined based on the noise level.

The second term is the new proposed DSP involved term aforementioned, and $K_n \varphi$ is the forward model of DPC which can be regarded as generating DPC images with a given estimation of $\varphi$ and $K_n$. The $L_0$-norm is used to achieve sparsity promotion. The third term is the Hessian regularization which can be regarded as higher order TV regularization. The Hessian regularization uses the second gradient which has smoother regularization result than that of TV one which uses the first order gradient. Here $\mathcal{H}$ denotes the Hessian gradient operator.

Since Eq. (9) involves $L_0$-norm which is a NP-hard problem, we develop two frameworks to approximately tackle the $L_0$-norm term.

*4.1 Half quadratic splitting (HQS) framework*

Equation (9) can be solved by the half quadratic splitting (HQS) algorithm in an alternating manner [23, 24]. By introducing $N + 1$ auxiliary variables, $\psi_n$, $(n = 1, 2, \cdots, N)$ and $G$ with respect to the $K_n\varphi$ and $\mathcal{H}\varphi$ respectively, the problem in Eq. (9) is converted into $N + 2$ sub-problems which are

$$\begin{cases} \arg\min_{\varphi} \sum_{n=1}^{N} \|K_n\varphi - s_{n,dpc}\|_2^2 + \alpha_0 \sum_{n=1}^{N} \|K_n\varphi - \psi_n\|_2^2 + \beta_0 \|\mathcal{H}\varphi - G\|_2^2 \\ \arg\min_{\psi_n} \alpha_0 \|K_n\varphi - \psi_n\|_2^2 + \alpha \|\psi_n\|_0, \quad n = 1, 2, \cdots N \\ \arg\min_{G} \beta_0 \|\mathcal{H}\varphi - G\|_2^2 + \beta \|G\|_1 \end{cases} \quad (10)$$

Here $\alpha_0$ and $\beta_0$ are sufficient large penalty parameters such that they enforce their corresponding L$_2$-norm terms approach to zero. For $\psi_n$ sub-problems (L$_0$-norm), it can be approximately solved by hard-threshold, and the closed-form solution is given by [24]

$$\psi_n = \begin{cases} K_n\varphi, & \sum_{n=1}^{N} |K_n\varphi|^2 > \alpha/\alpha_0 \\ 0, & else \end{cases} \quad (11)$$

For $G$ sub-problem (L$_1$-norm), it is solved by soft-threshold and the closed-form solution is given by [23]

$$G = \mathrm{sign}(\mathcal{H}\varphi) \circ \max(|\mathcal{H}\varphi| - \beta/\beta_0, 0). \quad (12)$$

$\circ$ denotes the element-wise multiplication. After $G$ and $\psi_n$ are obtained, the $\varphi$ sub-problem is pure-quadratic and can be solved by setting the derivative with respect to $\varphi$ to zero, and the closed-form solution is given by

$$\varphi = \frac{\sum_{n=1}^{N} K_n^T (s_{n,dpc} + \alpha_0 \psi_n) + \beta_0 \mathcal{H}^T G}{(1+\alpha_0) \sum_{n=1}^{N} K_n^T K_n + \beta_0 \mathcal{H}^T \mathcal{H}}, \quad (13)$$

where $T$ denotes the transpose of matrix. Since $K$ and $\mathcal{H}$ are block-circulant and can be digonalized by 2D Fourier transform, $\varphi$ can be computed directly and only at the cost of Fourier transform.

To make the algorithm practically working well, a trick is to let and $\alpha_0$ and $\beta_0$ increase in each iteration to gradually enforce $\|K_n\varphi - \psi_n\|_2^2 \to 0$ and $\|\mathcal{H}\varphi - G\|_2^2 \to 0$, and the $\varphi$, $\psi_n$ and $G$ are updated iteratively with known of each other. In such a manner, the steps of the half quadratic splitting algorithm are sketched in **Algorithm 1**. In an actual implementation $a = 2$ by default. $\alpha_{max} = 10^3$ and $\beta_{max} = 10^5$ to ensure that $\alpha_0$ and $\beta_0$ become large values.

This half quadratic splitting is one of the standard methods for solving regulated optimization problem with L$_0$-norm penalty and is widely used in many blind and non-blind deconvolution tasks. However, it is not fully optimized for DPC tasks as it involves large amounts of iterations before $\alpha_0$ and $\beta_0$ reach the $\alpha_{max}$ and $\beta_{max}$ since the initial values of $\alpha_0$ and $\beta_0$ are usually small [23]. For example, if the initial value of $\beta_0$ is 0.1, it takes 20 iterations to finish the loop. For $\alpha_0$ it takes 14 iterations. Therefore, it takes a total of 280 iterations for the final $\varphi$.

In order to boost the optimizing of Eq. (9), we propose, in next subsection, to solve Eq. (9) under a framework of Richardson-Lucy deconvolution.

---

**Algorithm 1**: Half quadratic splitting for Eq. (10)

**Input**: DPC image $s_{n,dpc}$, point spread function $K_n$, penalty parameters $\alpha$ and $\beta$
$\alpha_0 \leftarrow \alpha$, $\varphi = 1$
**While** $\alpha_0 < \alpha_{\max}$ **do**
    Solving $\psi_n$ sub-problem via Eq. (11) with given $\varphi$
    $\beta_0 \leftarrow \beta$
    **While** $\beta_0 < \beta_{\max}$ **do**
        Solving $G$ sub-problem via Eq. (12) with given $\varphi$
        Solving $\varphi$ sub-problem via Eq. (13) with given $\psi_n$ and $G$
        $\beta_0 \leftarrow a \cdot \beta_0$
    **End while**
    $\alpha_0 \leftarrow a \cdot \alpha_0$
**End while**
**Output**: Quantitative phase image $\varphi$

---

*4.2 Richardson-Lucy deconvolution (RLD) framework*

The data fidelity term in Eq. (9) is corresponding to the Richardson-Lucy deconvolution (RLD) under the Gaussian noise assumption [25]. We are not using the Poison noise assumption since it requires the norm of kernel being equal to 1, while in the DPC case the norm of kernel is 0 according to Eq. (6).

The Richardson-Lucy deconvolution uses the gradient descent to iteratively search for the optimal solution for deconvolution results. However, the L$_0$-norm term in Eq. (9) is non-differentiable which makes the gradient of $E(\varphi)$ impossible to be calculated.

An expedient measure to tackle with the L$_0$-norm is to replace it with L$_1$-norm as the L$_1$-norm is the convex approximation of L$_0$-norm. However, the power of L$_1$-norm on sparse regularization is not as strong as that of L$_0$-norm [26], which makes us searching for better candidate for substituting of L$_0$-norm, which is

$$f(x) = 1 - \exp(-c|x|), \quad \lim_{c \to \infty} f(x) \to \|x\|_0, \tag{14}$$

The evolution of $f(x)$ with respect to parameter $c$ is shown in **Supplementary Note 2** where this function shows its ability to approximate the function of L$_0$-norm and most importantly, it is differentiable (except for origin).

Therefore, we can choose $f(x)$ to approach the L$_0$-norm by assign an appropriate value of $c$ ($c = 10$ in this study), and back propagate its derivative to update the cost function. Replacing the L$_0$-norm term in Eq. (9) by Eq. (14) yielding

$$E(\varphi) = \sum_{n=1}^{N} \|K_n \varphi - s_{n,dpc}\|_2^2 + \alpha \sum_{n=1}^{N} \left[1 - \exp(-c|K_n \varphi|)\right] + \beta \|\mathcal{H}\varphi\|_1. \tag{15}$$

Then the derivative (sub-gradient) of $E$ with respect to $\varphi$ is given by

$$\frac{\partial E}{\partial \boldsymbol{\varphi}} = 2\sum_{n=1}^{N} \boldsymbol{K}_n^T \left( \boldsymbol{K}_n \boldsymbol{\varphi} - \boldsymbol{s}_{n,dpc} \right) + \alpha c \sum_{n=1}^{N} \boldsymbol{K}_n^T \left[ \exp\left(-c|\boldsymbol{K}_n \boldsymbol{\varphi}|\right) \frac{\boldsymbol{K}_n \boldsymbol{\varphi}}{|\boldsymbol{K}_n \boldsymbol{\varphi}|} \right] + \beta \mathcal{H}^T \left( \frac{\mathcal{H} \boldsymbol{\varphi}}{|\mathcal{H} \boldsymbol{\varphi}|} \right), \quad (16)$$

In a conventional gradient descent framework, $\boldsymbol{\varphi}$ at $t$-th iteration is updated using gradient descent $\boldsymbol{\varphi}_{t+1} = \boldsymbol{\varphi}_t - \eta \cdot \partial_\varphi E$, and $\eta$ is the step size. However, gradient bursting problem happens for conventional gradient descent method since the absolute value of gradient of $f(x)$ increases volatilely for small value of $\boldsymbol{K}_n \boldsymbol{\varphi}$. With an extreme large gradient value, the gradient descent method easily misses an optimal point even if a small step size is given.

To avoid this gradient bursting problem and ensure the optimizing accuracy and speed, we adopt the adaptive momentum estimation with Nesterov acceleration (N-Adam) [27] and $\boldsymbol{\varphi}$ is updated according to

$$\begin{cases} \boldsymbol{v}_t = \rho_1 \cdot \boldsymbol{v}_{t-1} + (1-\rho_1) \cdot \frac{\partial L}{\partial \boldsymbol{\varphi}} \\ \boldsymbol{\upsilon}_t = \rho_2 \cdot \boldsymbol{\upsilon}_{t-1} + (1-\rho_2) \cdot \left(\frac{\partial L}{\partial \boldsymbol{\varphi}}\right)^2 \\ \hat{\boldsymbol{v}} = \boldsymbol{v}_t / (1-\rho_1^t) \\ \hat{\boldsymbol{\upsilon}} = \boldsymbol{\upsilon}_t / (1-\rho_2^t) \\ \boldsymbol{\varphi}_{t+1} = \boldsymbol{\varphi}_t - \eta \cdot \frac{1}{\sqrt{\hat{\boldsymbol{\upsilon}}} + \xi} \left[ \rho_1 \cdot \hat{\boldsymbol{v}} + (1-\rho_1) \cdot \frac{\partial L}{\partial \boldsymbol{\varphi}} \right] \end{cases}, \quad (17)$$

where $\boldsymbol{v}_t$ ($\boldsymbol{v}_t = \boldsymbol{0}$) is the first momentum estimation, and $\boldsymbol{\upsilon}_t$ ($\boldsymbol{\upsilon}_t = \boldsymbol{0}$) is the second-momentum estimation. $\hat{\boldsymbol{v}}$ and $\hat{\boldsymbol{\upsilon}}$ are bias-corrected of $\boldsymbol{v}_t$ and $\boldsymbol{\upsilon}_t$ respectively. $\rho_1$ and $\rho_2$ are decay rates for the momentum estimations. $\rho_1 = 0.9$ and $\rho_2 = 0.999$ are used in the following experiments. $\eta = 0.05$ is the learning rate. $\xi$ is a small value used to avoid the condition of dividing by 0.

---

**Algorithm 2**: N-Adam Richardson-Lucy deconvolution for Eq. (15)

**Input**: DPC image $s_{n,dpc}$, $\boldsymbol{K}_n$, initialize $\boldsymbol{\varphi}_1$, $\boldsymbol{v}_0 = \boldsymbol{\upsilon}_0 = 0$, $t = 1$.
**While** $t < t_{max}$ **do**
    Calculate gradient $\partial_\varphi E$ using Eq. (16)
    Updating $\boldsymbol{v}_t$, $\boldsymbol{\upsilon}_t$, and $\boldsymbol{\varphi}_{t+1}$ using Eq. (17)
    $t \leftarrow t+1$
**End while**
**Output**: Quantitative phase image $\boldsymbol{\varphi}$

---

The original Adam allows the momentum gradually accumulated and won't go extrema fast even if a large gradient appears in our approximated $L_0$-norm optimization. Combined with Nesterov acceleration, the N-Adam is able to converge faster and more stable than that the original Adam algorithm.

In such a manner, the steps of the Richardson-Lucy deconvolution with N-Adam gradient are sketched in **Algorithm 2**. $t_{max} = 150$ for most of our experiments, and more iteration may be required for better results. The initial guess of $\varphi_1$ can be calculated using L$_2$-norm regularization with $\alpha = 1$.

*4.3 Adaptive noise sensor*

In our model, there are two penalty parameters $\alpha$ and $\beta$ that need to be manually adjusted, which limits our model implementation in various conditions. In principle, $\alpha$ and $\beta$ are adjusted according to the noise level of the capture images. Inspired by this, we designed an adaptive noise sensor to evaluate the noise levels of the images [28]. This sensor is depicted by

$$\Lambda = \frac{1}{20}\sqrt{\frac{\pi}{2}}\frac{1}{N}\frac{1}{WH}\sum_{n=1}^{N}\sum_{x=1}^{W}\sum_{y=1}^{H}\left|s_n(x,y)\otimes\mathcal{L}\right|, \tag{18}$$

where $\mathcal{L} = [-1, 2, -1; -2, 4, -2; -1, 2, -1]$ is the Laplacian operator.

For noise-free images, the Laplacian operator is able to extract the edge information. On the contrary, for noise corrupted images, Laplacian operator is able to enlarge the noise signals as it is sensitive to noise pixels, and the average absolute pixel value of the filtered image can be regarded as a measurement of noise level.

Accordingly, the parameters are given by

$$\begin{cases} \alpha = \Lambda \\ \beta = \alpha/2 \end{cases} \tag{19}$$

In some cases, parameter $\beta$ may need manual adjustment to achieve better deconvolution results.

## 5. Experimental results

*5.1 Validation with simulated phase sample*

First we validated our algorithm on simulated data. System parameters were $\lambda = 0.530$ μm, $NA = 0.3$, a magnification of $\times 10$ and pixel size of the camera of 3.46 μm. A binary phantom object for the phase target is shown in Fig. 3 (a) where the phase is ranged in $[0, 2]$ rad.

According to the forward model, we generated 2 DPC-images corresponding to top-bottom illuminations and left-right illuminations and added Gaussian noises on them. The standard deviation of the Gaussian noise is $0.2(I_{max} - I_{min})$, where $I_{max}$ and $I_{min}$ are maximum and minimum value of each DPC image. Figure 3 (b) shows one of the DPC image where the image is severely corrupted by noise. The performance of our SH-DPC is compared against the following algorithms:
1) Tikhonov (L$_2$-norm) regularization based DPC (L$_2$-DPC): we use the traditional DPC reconstruction scheme by setting $\alpha = 0.0001$ in Eq. (1) with L$_2$-norm penalty.
2) TV-regularization based DPC (TV-DPC): we use the gradient L$_1$-norm penalty function, and $\alpha = \Lambda$ according to Eq. (15).

The reconstruction phase for L$_2$-norm regularization is shown in Fig. 3 (c). We also plot the phase value along the red curve in Fig. 3 (g) to shown detailed comparison between reconstruction results and ground-truth. According to Fig. 3(c) and Fig. 3 (g), the noise signal is still very large and the algorithm fails to reconstruct the correct phase image.

Results of TV-DPC are shown in Fig. 3(d) and Fig. 3 (h). The penalty parameter $\alpha$ are given by Eq. (15), where $\Lambda = 0.123$. The TV-regularization suppresses the impact of the noise signal as the reconstructed result approaches to the ground truth as shown in Fig. 3 (h). However, the phase image is still corrupted by noise signals as the 'white-fog effects' present in the background.

Results of proposed SH-DPC using HQS solution are shown in Fig. 3(e) and Fig. 3 (i). According to Fig. 3 (e), the noise signal is significantly suppressed as there is no 'white-fog effects' in the background. Meanwhile, data fidelity is also guaranteed according to Fig. 3 (i). As for RLD solution, our proposed $L_0$-approximation successfully achieves the sparse penalty, as the 'white-fog effects' in the background are removed as shown in Fig. 3(f). The reconstruction phase is consistent with ground-truth as shown in Fig. 3 (j).

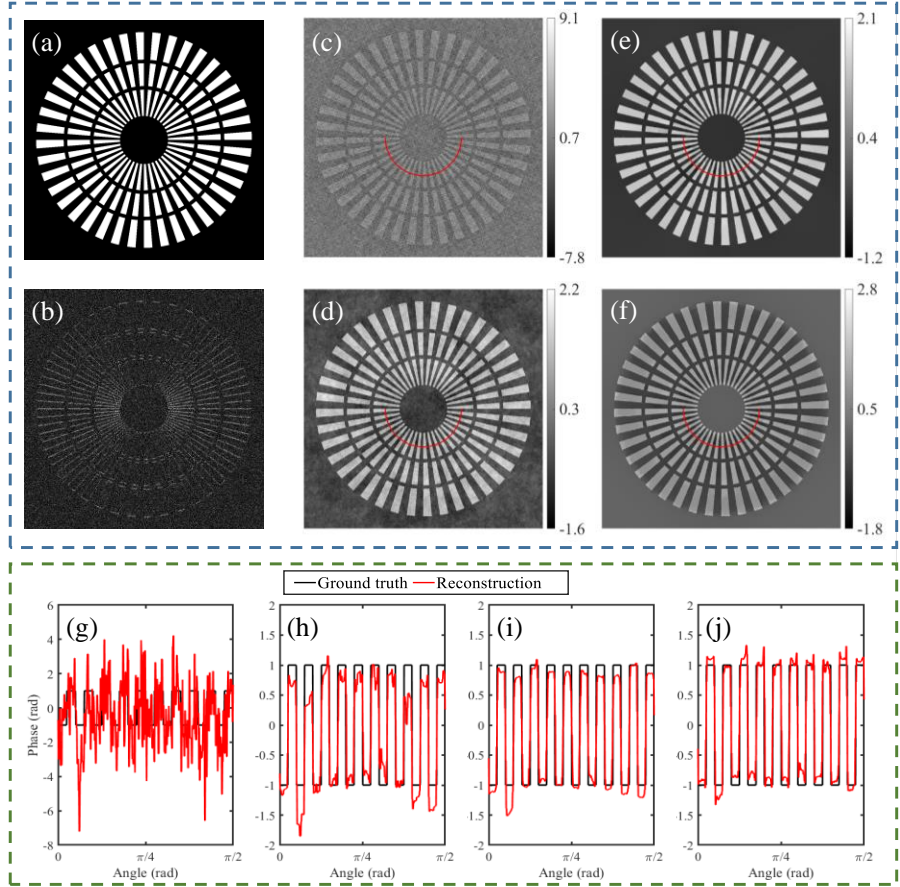

Fig. 3, simulation of DPC reconstruction. (a) ground-truth of the phase phantom. (b) Absolute value of DPC image corrupted by Gaussian noise. (c) Reconstruction results of $L_2$-DPC. (d) Results of TV-DPC. (e) Results of proposed model using HQS. (f) Results of proposed model using RLD. Quantitative phase plot along the red arc. (g) $L_2$-norm regularization. (h) TV regularization. (i) Proposed model using HQS. (j) Results of proposed model using RLD.

We then performed 5 groups of simulation experiment using the phase maps shown in Tab. 2. All the ground-truth phase maps are of size (600 × 600) pixels, and have values in range of [0, 1] radians. For each ground-truth image, we generated the 2 DPC-images corresponding to top-bottom illuminations and left-right illuminations according to the forward model of DPC, and corrupted the DPC images using Gaussian noise with assigned signal-to-noise-ratio (SNR).

We quantitatively measured the quality of reconstructed phase maps by computing the linear regressed SNR (LSNR) to get rid of additive constants $b$. Between the ground-truth map $\varphi$ and a reconstructed one $\tilde{\varphi}$, the LSNR measure is defined as

$$\text{LSNR}(\varphi, \tilde{\varphi}) = \max_{b \in \mathbb{R}} \quad 10\log\left[\frac{\|\varphi\|_2^2}{\|\varphi - (\tilde{\varphi} + b)\|_2^2}\right]. \tag{20}$$

According to Eq. (20), larger value of LSNR implies better phase reconstruction quality. For each level of SNR, we further duplicate the experiment by 10 times, and calculate the average LSNR for a reliable comparison. Results are listed in Tab. 1.

As listed in Tab. 1, the SH-DPC method with two solutions outperforms the other algorithm in almost all of the cases. Especially for moderate (SNR = 15, 20) and high levels of noises (SNR = 10), the SH-DPC significantly improves the phase reconstruction quality.

Table 1. Phase reconstruction quality of the algorithms

| Test Images (5 groups) | SNR | Method | | | |
|---|---|---|---|---|---|
| | | $L_2$-DPC | TV-DPC | DSP-HQS | DPS-RLD |
| 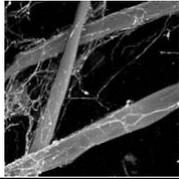 | 10 | -29.156 ± 0.174 | 24.438 ± 0.481 | 30.224 ± 0.484 | 25.903 ± 0.329 |
| | 15 | -17.804 ± 0.101 | 28.895 ± 0.428 | 33.457 ± 0.388 | 29.788 ± 0.283 |
| | 20 | -6.300 ± 0.165 | 32.492 ± 0.381 | 35.841 ± 0.272 | 33.074 ± 0.229 |
| | 30 | 15.973 ± 0.102 | 37.360 ± 0.176 | 38.972 ± 0.115 | 37.582 ± 0.120 |
| 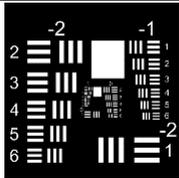 | 10 | -28.092 ± 0.173 | 29.843 ± 0.310 | 35.600 ± 0.195 | 34.164 ± 0.176 |
| | 15 | -16.513 ± 0.131 | 33.826 ± 0.167 | 37.242 ± 0.149 | 36.028 ± 0.114 |
| | 20 | -5.192 ± 0.173 | 36.136 ± 0.258 | 37.627 ± 0.089 | 36.989 ± 0.093 |
| | 30 | 16.788 ± 0.111 | 38.045 ± 0.090 | 37.937 ± 0.044 | 37.956 ± 0.047 |
| 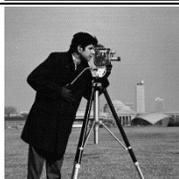 | 10 | -18.409 ± 0.115 | 36.665 ± 0.835 | 43.552 ± 0.626 | 40.934 ± 0.570 |
| | 15 | -6.954 ± 0.184 | 41.416 ± 0.645 | 45.929 ± 0.273 | 44.051 ± 0.243 |
| | 20 | 4.556 ± 0.117 | 44.334 ± 0.514 | 47.340 ± 0.378 | 46.149 ± 0.351 |
| | 30 | 26.491 ± 0.150 | 47.543 ± 0.192 | 47.841 ± 0.090 | 48.010 ± 0.104 |
| 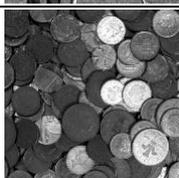 | 10 | -21.047 ± 0.142 | 33.136 ± 0.608 | 36.784 ± 0.670 | 33.749 ± 0.493 |
| | 15 | -9.666 ± 0.165 | 37.355 ± 0.304 | 39.981 ± 0.256 | 37.636 ± 0.224 |
| | 20 | 1.882 ± 0.214 | 40.421 ± 0.323 | 42.262 ± 0.252 | 40.614 ± 0.226 |
| | 30 | 23.923 ± 0.129 | 43.660 ± 0.147 | 44.155 ± 0.163 | 43.782 ± 0.146 |
| 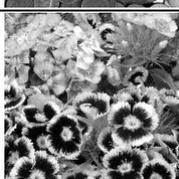 | 10 | -18.021 ± 0.142 | 32.370 ± 0.374 | 36.521 ± 0.575 | 35.806 ± 0.349 |
| | 15 | -6.511 ± 0.174 | 36.693 ± 0.379 | 39.703 ± 0.423 | 38.892 ± 0.294 |
| | 20 | 4.869 ± 0.121 | 40.275 ± 0.208 | 42.033 ± 0.219 | 41.418 ± 0.171 |
| | 30 | 26.738 ± 0.165 | 44.755 ± 0.087 | 45.344 ± 0.093 | 45.310 ± 0.089 |
| Average (100 images) | 10 | -22.945 ± 4.815 | 31.290 ± 4.138 | 36.536 ± 4.317 | 34.111 ± 4.885 |
| | 15 | -11.490 ± 4.821 | 35.637 ± 4.213 | 39.262 ± 4.125 | 37.279 ± 4.608 |
| | 20 | -0.037 ± 4.840 | 38.732 ± 4.112 | 41.021 ± 4.070 | 39.649 ± 4.356 |
| | 30 | 21.983 ± 4.736 | 42.272 ± 3.989 | 42.850 ± 3.835 | 42.528 ± 4.041 |

Moreover, the HQS method gains higher LSNR score than the RLD. The main reason is that HQS is using hard-threshold method to approximate the $L_0$-norm. The hard-threshold directly assigns zeros to the matrix where the pixel value is less than the threshold. While the RLD is using gradient descent, during this process, the value of the pixels can extremely approach to zeros but cannot be equal to zero. Therefore, the sparsity of the results of RLD method is less than that of HQS method. On the contrary, the RLD requires less iterations which is 2-times more efficiency than the HQS method.

*5.2 Experimental verification: quantitative phase target*

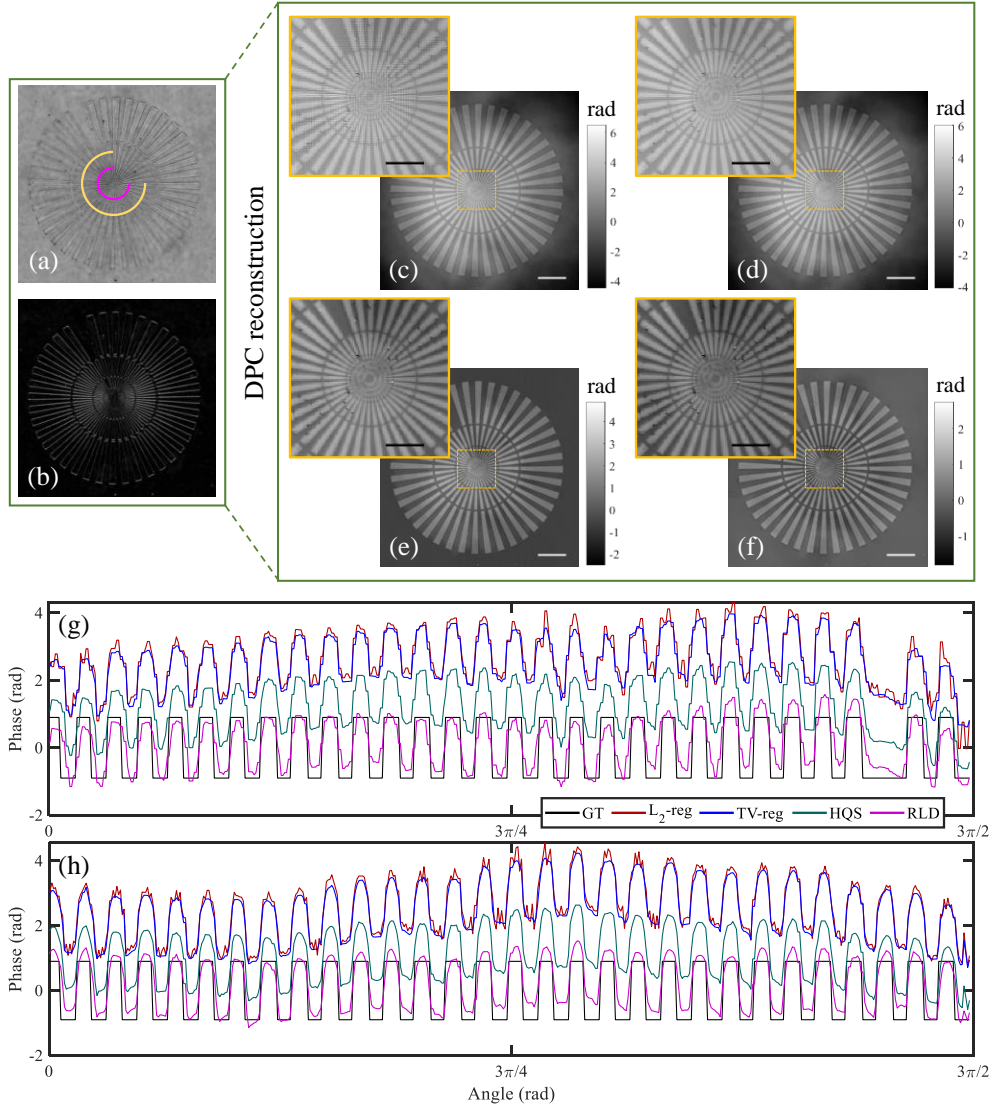

Fig. 4. Reconstruction results for focal-star pattern (a1) $L_2$-norm regularization. (b1) TV regularization. (c1) Proposed model using HQS. (d1) Results of proposed model using RLD. The scale bar is 65 μm. (a2) to (d2) are zoomed in part for the region labelled by yellow box in (a1). The scale bar is 20 μm. (e) and (f) are plot of phase results along (e) the magenta arc in Fig. 4(a1)-(d1). (f) Along the yellow arc in Fig. 4(a1)-(d1).

We further perform DPC experiment to show the experimental robustness of our proposed SH-DPC. First, we validate our SH-DPC among using quantitative phase target (QPT, BenchMark, USA). The focal star is 100nm height, and the refractive index is 1.52. The QPT is imaged by an objective lens with $NA = 0.3$ and magnification $\times 10$. The pixel size of the camera (Edge 4.2 mono, PCO, German) is 6.5 μm. We image the QPT using a ring shaped LED array (WS2812B, China), and the illumination pupil is also a ring shape whose outer radius is $NA/\lambda$, and inner radius is $0.9\, NA/\lambda$.

We collected 4 oblique illumination images corresponding to top, bottom, right, and left illumination, therefore we obtain 2 groups of DPC images which is used for phase reconstruction. Our reconstruction results are compared against with the results constructed using $L_2$-regularization and TV-regularization.

As shown in Figs. 4(c), the phase is reconstructed using $L_2$-regularization where $\alpha = 10^{-5}$, and the signal is corrupted by noises spots (produced by deconvolution of noise pixels) and "white-cloud" effect (caused by unevenly illumination pattern generated by different groups of LEDs dot). With TV-regularization as shown in Fig. 4(d), the noises spots is suppressed but the "white-cloud" effect is still obvious.

SH-DPC reconstruction results are shown in Fig. 4(e) for HQS method, and 4(f) for RLD method. The "white-cloud" effect is effectively removed as the background component is uniformly distributed. The noise spots are also suppressed. By comparing the zoomed in part in for the small structures in Fig. 5(c2) and 5(d2), our proposed SH-DPC algorithm will neither cause the loose of information nor will create unexpected artefacts.

We plot the quantitative phase profile along the yellow and magenta lines shown in Figs. 4(e) and 4(f), respectively. As shown in Fig. 4(e) and 4(f), results for $L_2$-DPC (red curves) and TV-DPC (blue curves) are divergent from the ground truth value (black curves) due to the impact of noise and 'white-cloud' background. The blue curves is smoother than the red curves implying the noise suppression effect of TV regularization. The SH-DPC reconstructions fits well with the ground truth as the green and magenta curves are close to the black curve in both Figs. 4(e) and 4(f), since the "white-cloud" background is removed. More experimental results are found in **Supplementary Note 3**.

*5.3 Experimental verification: Microsphere clusters*

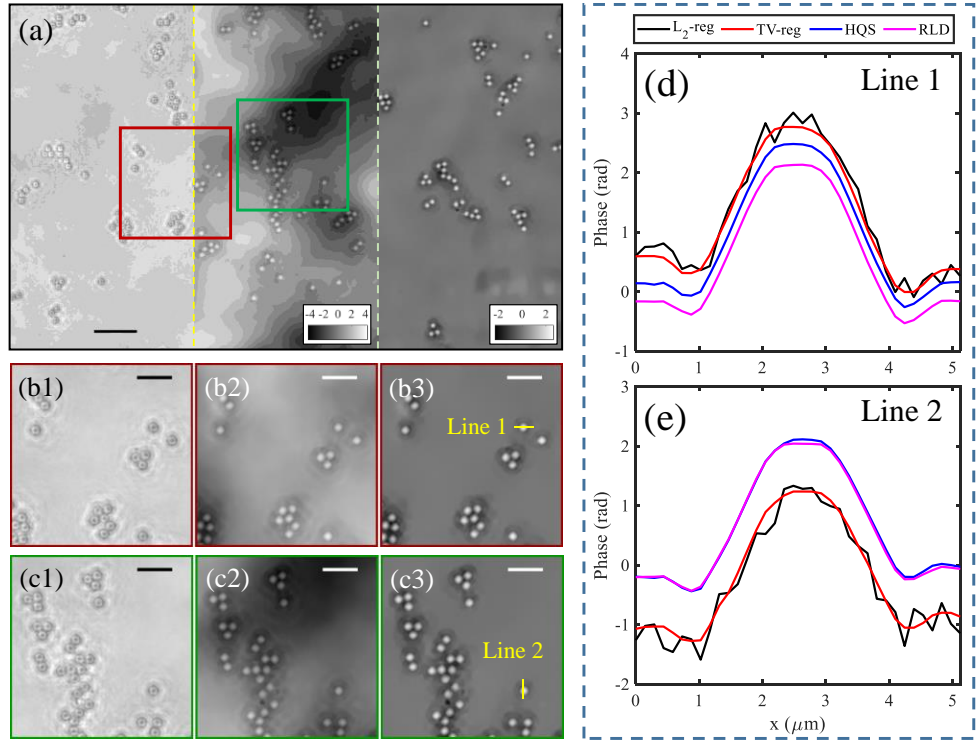

Fig. 5. Reconstruction results for microspheres cluster. (a) Montage of bright field image, TV-DPC result, and DSP-HQS results, Scale bar is 24 μm. (b1) to (b3) are enlarged parts for region in the green box in (a), and reconstructed by TV-DPC, HQS, and RLD methods. Scalar bar is 12 μm (c1) to (c3) are enlarged parts for region in the blue box. (d) and (e) are quantitative phase profile along Line 1 and Line 2, respectively.

We further validate our SH-DPC on the microspheres (Duke, USA) cluster sample. The microsphere is of 3um in diameters. The refractive index of the microsphere is 1.59 and is merged in cedar wood oil (RI=1.51). The sample is imaged with an objective lens of $NA = 0.75$, $\times 40$. The pixel size of the camera (SUA231GM-T, MindVision, China) is 5.85 μm.

The montage of full-field TV-DPC reconstruction (top-right part) and HQS reconstruction (bottom-right left) is shown in Fig. 5(a). The reconstruction by TV-DPC is severely degenerated by the "white cloud" effect due to the uneven illumination and noise signal among the full field of view. The SH-DPC removes the white cloud and noise effect as the background is clear and uniformly distributed. Zoomed in pictures in Fig. 5(b1) to 5(c3) shows the detailed comparison between bright field image, TV-DPC, and SH-DPC with HQS solution. Results from $L_2$-DPC and SH-DPC with RLD method are not shown.

Figure 5(d) plots the quantitative phase profiles along line1 in Fig. (b3), where the phase difference between peak and valley for four methods are consistent with each other. The curve for $L_2$-DPC is corrupted by noise signal and is not as smooth as that of TV-DPC and SH-DPCs'. Similar analyses applied to Fig. 5(e) for the phase profile along line 2 in Fig. 5(c3). More experimental results are found in **Supplementary Note 4**.

*5.4 Validation*

We first validate our SH-DPC by imaging the unstained gastric cancer cells. The cells are imaged by an objective with $NA = 0.3$, $\times 10$. The pixel size of the camera is 6.5 μm. The illumination wavelength is 532 nm. A region of interest (ROI) of size (1200 × 1200) is chosen as shown in Fig. 6 (a). The bright-field image has very low contrast, as the gastric cancer cells can be regarded as pure-phase objects.

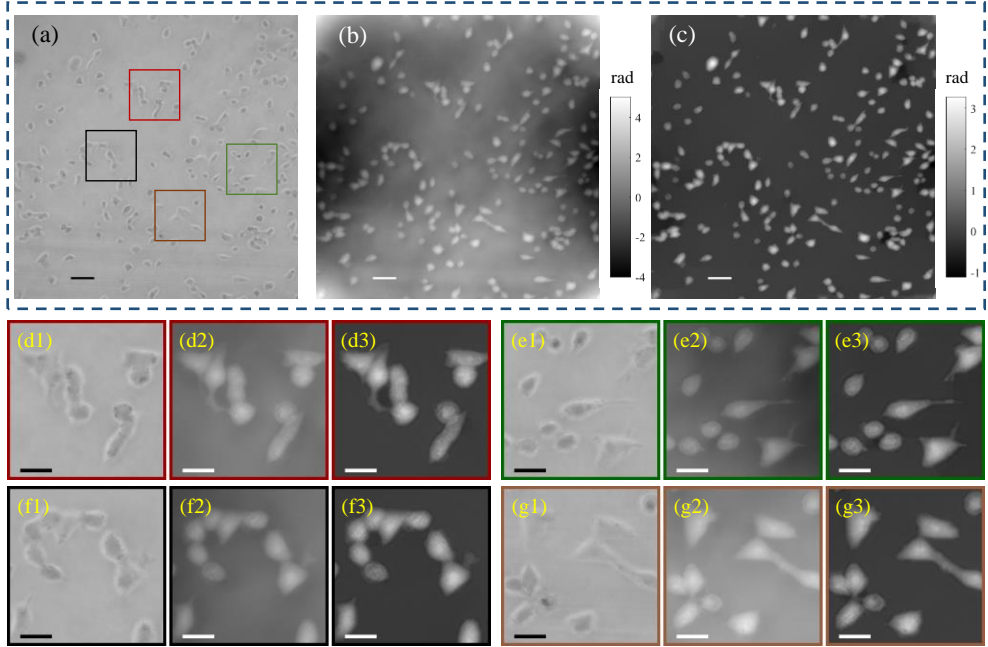

Fig. 6. Reconstruction results for gastric cancer cells. *β* = *α*. (a) Bright field image. (b) and (c) are recovery results using TV-reg and our proposed method using RLD. Scale bar is 65 μm. (d1) to (g3) are zoomed-in images corresponding to the areas in the green, black, orange, and red boxes. Scale bar is 30 μm

Based on previous results, we compared the SH-DPC with RLD algorithm against TV-DPC. As shown in Fig. 6 (b), the TV-DPC cannot suppress the 'white-fog' effect caused by uneven illumination fluctuation, while our proposed SH-DPC significantly removes the 'white-fog' effect for entire ROI as shown in Fig. 6 (c). Zoomed-in images listed in Figs. 6 (d1) to 6 (g3) show that our SH-DPC improves the reconstruction quality while maintain the data fidelity, as the small structures such as the edges and branches of the cells are preserved.

We also tested SH-DPC on color multiplexed single-shot DPC (cDPC) as shown in Fig. 7 (a). The samples (HeLa cells) were illuminated by red, blue, and green colored LED and were imaged by an objective lens with $NA = 0.3$, $\times 10$. The pixel size of the camera is 5.85 μm. Figure 7 (a) is the input raw image for DPC reconstruction. Since it is illuminated by oblique plane wave, the phase contrast effect can be observed.

Since the DSP is a universal prior for DPC image, it also works for cDPC as the 'white-fog' effect is suppressed shown in Fig. 7 (c). Zoomed-in images listed in Fig. 7 (d1) to 7 (g3) show that our SH-DPC improves the reconstruction quality where the clarity of the image is increased, and more delicate structures can be overserved.

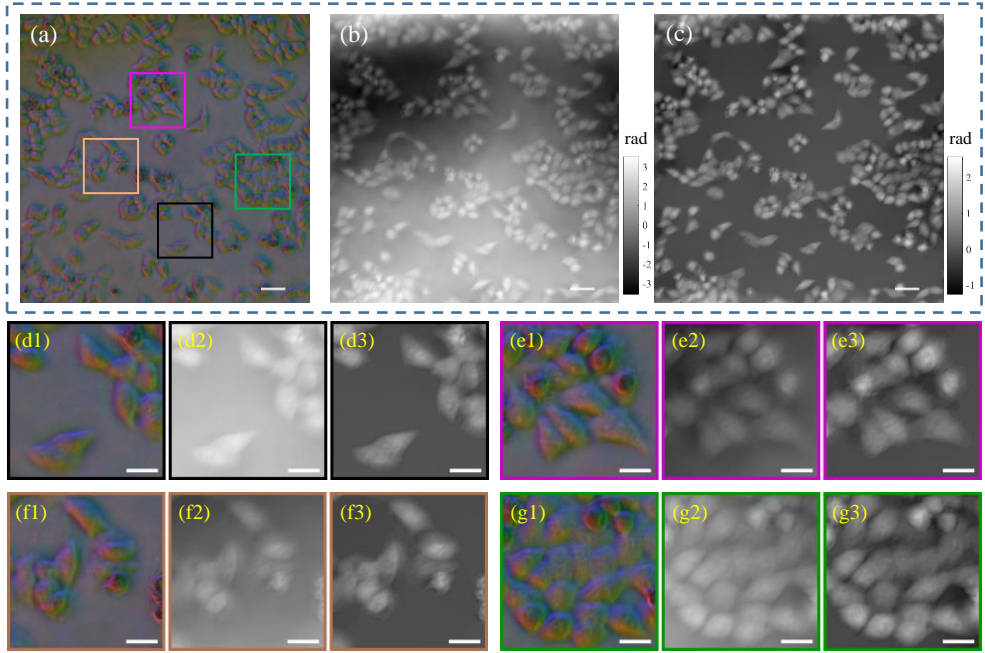

Fig. 7. Reconstruction results for HeLa cells using cDPC. *β* = *α*. (a) Raw image for DPC reconstruction. (b) and (c) are recovery results using TV-DPC and our proposed method using RLD. Scale bar is 65 μm. (d1) to (g3) are zoomed-in images corresponding to the areas in the black, purple, orange, and green boxes. Scale bar is 30 μm.

Final, we show that our proposed SH-DPC benefits some potential applications using phase imaging, such as cell segmentation tasks (Automatic delineation of the cell boundaries). We applied Wang's method [26, 29] on both Fig. 7 (b) and Fig. 7 (c) with identical model parameters, and segmentation results are shown in Fig. 8. As shown in Figs. 8 (a1) to 8 (c1), the segmentation results on TV-DPC output are not satisfactory due to lack of image contrast since the boundary signal of cells are reduced by the haze effect. Only few cells are labeled as shown in Fig. 8 (c1). On the contrary, the segmentation results on DSP-DCP are rather satisfactory where the boundary of the cells are recognized as shown in Fig. 8 (a2), and most of the cells are labeled as shown in Fig. 8 (c2).

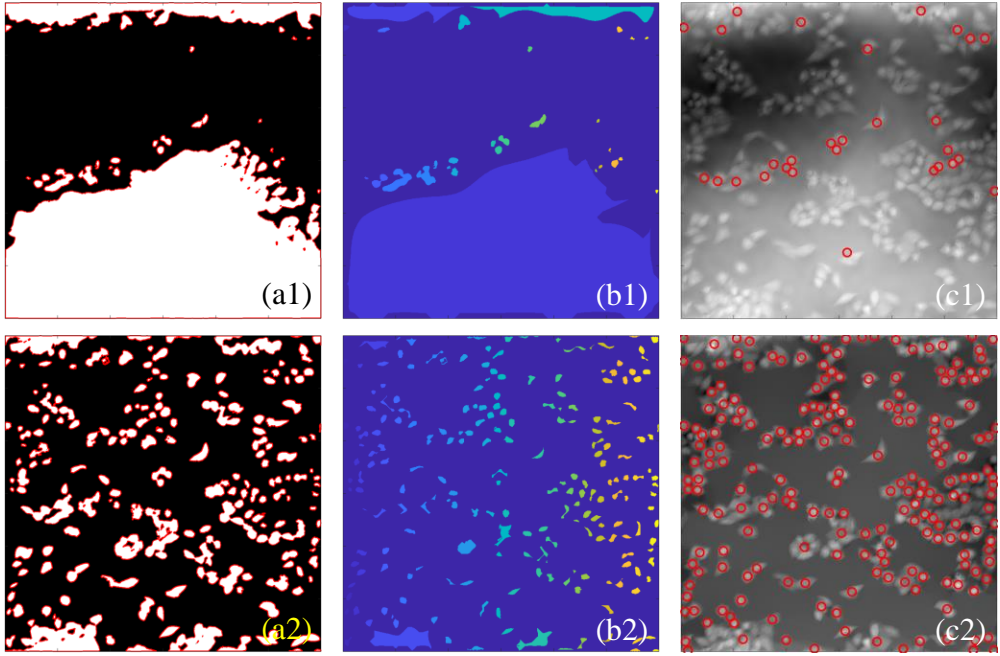

Fig. 8. Cell segmentation using phase image. (a1), (b1) and (c1) are boundary map, segmentation map, and counting map for image in Fig. 7 (b) (TV-reg). (a2), (b2) and (c2) are boundary map, segmentation map, and counting map for image in Fig. 7 (c) (SH-DPC).

## 6. Analysis and discussions

It is not surprising that the dark-field sparse prior (DSP) enables us to design a DPC deconvolution algorithm that outperforms state-of-the-art methods, but also obtains competitive results on phase recovery quality even for images corrupted by heavy noises and un-even background intensity. The DSP significantly increases the experimental robustness without the need of modification of imaging system. In this section we further analyze the DSP, compare it with related method, and discuss its potential applications.

### 6.1 Effectiveness of the dark-field sparse prior

Our deconvolution model without the DSP reduces to the deconvolution with high-order total variation regularization (Hessian regularization). To ensure fair comparison, we disable DSP in our implementation. The input DPC image are shown in Fig. 9 (a) and 9 (b). The DPC image are taken log10 to clearly show the noise background due to the noises and uneven illuminations. Fig. 9 (d) and 9 (e) shows the DPC deconvolution using RLD method with only Hessian regularization, and DSP, respectively. As shown in Figs. 9 (e) and 9 (f), using DSP generates intermediate of $K_n \varphi$ with sparser background, which favor an idea DPC image, suppresses the noise signal, and facilitates the deconvolution process. Also, the latent $K_n \varphi$ becomes sparser with more iterations. Both of our simulation and experimental results concretely demonstrate the effectiveness of the DSP.

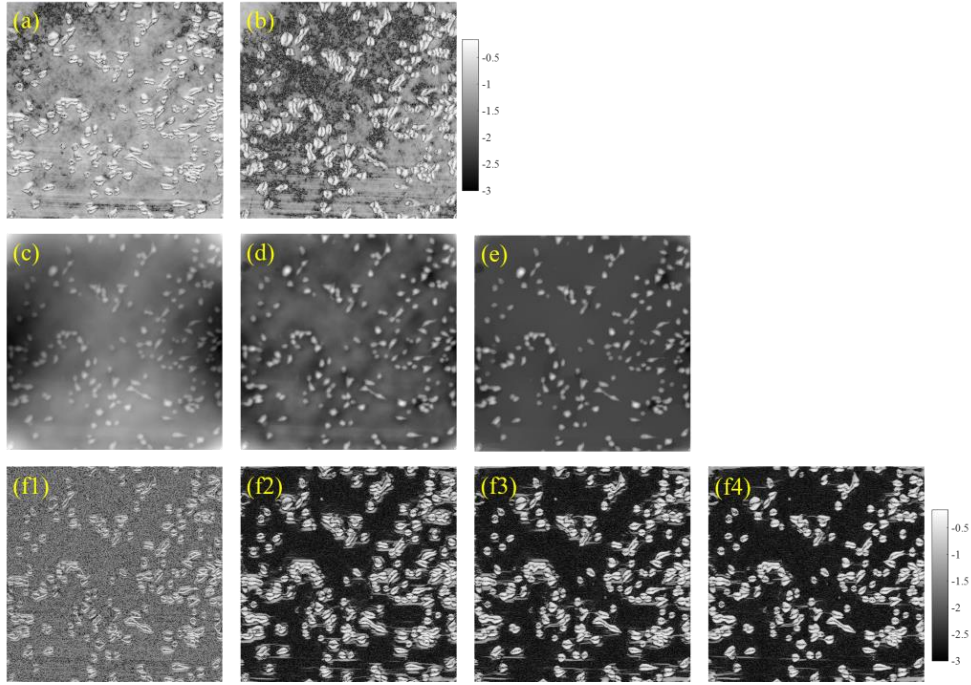

Fig. 9. DPC deconvolution. (a) and (b) are patterns for $\lg(|s_{1,dpc}|)$ and $\lg(|s_{2,dpc}|)$. Pixel value that are less than $10^{-3}$ are limited to $10^{-3}$. (c), (d) and (e) are phase reconstruction results using TV-reg, only Hessian-reg, and only DSP with RLD method, respectively. (f1) to (f4) are latent images for $\lg(|K_1\varphi|)$ at 1st, 25-th, 50-th, and 100-th iterations. Pixel value that are less than $10^{-3}$ are limited to $10^{-3}$.

### 6.2 Favored minimum of the cost function

The DSP is effective because it has lower energy for idea DPC images than for noise-corrupted DPC images. The DSP favors the idea DPC images. As shown in Figs. 9 (f), the Fig. 2(d), the idea DPC image is so sparser than the noise-corrupted DPC image that the DPC can simply distinguish the idea DPC image from that of noised DPC image. In this manner, the DSP gives lower energy to idea DPC images, and rectify the DPC deconvolution. It is worth noting that our proposed dark-field spares prior can be extended beyond the DPC since the DSP is also a universal feature for optical images as long as dark-field measurements are involved, for example, the Fourier ptychography [30]. Moreover, our proposed RLD method uses gradient descent to tackle the $L_0$-norm problem which is able to solve the non-convex cost function. These features make our DSP have great potential applications in computational imaging including but not limited to Fourier ptychography, diffraction tomography [31, 32], and fluorescence imaging [33].

### 6.3 Half-quadratic splitting vs Richardson-Lucy deconvolution

In this research, we proposed two frameworks to tackle the $L_0$-norm problem (1) the Half-quadratic splitting (HQS) and Richardson-Lucy deconvolution (RLD). It is worth noting that both frameworks provide approximate solutions to the $L_0$-norm. As the original $L_0$-norm is an NP-hard problem. The HQS uses the hard threshold to achieve the sparse promotion and has a stronger regularization ability on the sparsity than that of RLD method. However, the hard threshold may cause some sudden break-off on the cellular structures, which makes the results on real DPC data looks unnatural. Moreover, the HQS implies the optimization to be convex, which may not hold for real DPC experiment since many issues including illumination problems and vignetting effects contains non-linear processes which make the problem non-

**17** / 19

convex. While the RLD uses the N-Adam to update the cost function, it is more robust on the non-convex optimization and favors experimental conditions. The RLD is more efficient than HQS method as it requires fewer iterations. Empirically, we would recommend using RLD for DPC deconvolution.

## 7. Concluding remarks

In this research, we proposed a new image prior – the dark-field sparse prior (DSP), and formulated a non-convex model to facilitate the quantitative differential phase contrast inverse problem. We use $L_0$-norm to represent the DSP, and proposed two approximate solutions, (1) the HQS method and (2) the RLD method, to solve the NP-hard $L_0$-norm problem. Our adaptive parameter determination simplified parameter selection for the model, and both simulation and experimental results show the superiority of our proposed SH-DPC in terms of phase reconstruction quality and implementation efficiency as the SH-DPC is purely mathematical, which doesn't require any modification to the optical system compared to existing methods and benefit all DPC experiments. Finally, we illustrated that the proposed method is beneficial to image analysis such as cell segmentation.

**Funding:** This research was funded by the Key Research and Development Program of Anhui Province in China (grant No. 2022a05020028), the Natural Science Foundation of Anhui Province in China (grant No. 1908085MA14) and the project of Research Fund of Anhui Institute of Translational Medicine (grant No. 2021zhyx-B16).